\documentclass[reqno,11pt,twoside]{article}
\usepackage[hmargin=1.25in, vmargin=1.25in, a4paper, centering]{geometry}
\usepackage{amsmath, amsthm,amssymb}
\usepackage{fourier}
\usepackage{ebgaramond}
\usepackage{graphicx}
\graphicspath{{./}{../Figures-and-Scripts/}}
\usepackage{xcolor}
\usepackage[ font=small, labelfont=bf, labelsep=colon]{caption}
\usepackage{etoolbox}
\usepackage{longtable,booktabs,array}
\usepackage{algorithm}
\usepackage{algpseudocode}
\usepackage{marvosym}
\usepackage[section]{placeins}
\usepackage{tocloft}

\numberwithin{equation}{section}
\theoremstyle{definition}
\newtheorem{defi}{Definition}[section]

\theoremstyle{plain}
\newtheorem{theorem}[defi]{Theorem}
 \newtheorem{prop}[defi]{Proposition}
\newtheorem{lemma}[defi]{Lemma}

\theoremstyle{remark}
\newtheorem{remark}[defi]{Remark}
\newcommand{\addQEDstyle}[2]{\AtBeginEnvironment{#1}{\pushQED{\qed}\renewcommand{\qedsymbol}{#2}}\AtEndEnvironment{#1}{\popQED}}
\addQEDstyle{remark}{$\blacktriangleleft$}
\newcommand{\A}{\mathcal{A}}

\newcommand{\area}{\operatorname{Area}}

\newcommand{\Dir}{{\mathrm{D}}}
\newcommand{\Neu}{{\mathrm{N}}}
\newcommand{\dr}{\mathrm{d}}
\newcommand{\ir}{\mathrm{i}}
\newcommand{\er}{\mathrm{e}}

\renewcommand{\epsilon}{\varepsilon}

\newcommand{\disk}{\mathbb{D}}

\newcommand{\Lfloor}{\left\lfloor}
\newcommand{\Rfloor}{\right\rfloor}
\newcommand{\entire}[1]{\Lfloor #1 \Rfloor}
\newcommand{\ceiling}[1]{\left\lceil#1\right\rceil}

\renewcommand\footnotemark{}
\usepackage[nobottomtitles,pagestyles]{titlesec}
\titleformat{\section}
{\normalfont\large\bfseries}
{\filcenter\S\thesection.}{1ex}{\filcenter}
\widenhead*{0.5in}{0.5in}
\renewpagestyle{plain}[\bfseries]{\footrule\setfoot{}{\thepage}{}}
\newpagestyle{mystyle}[\bfseries]{
\headrule
\sethead[\thepage][][Aharonov--Bohm eigenvalues of the disk]{N. Filonov, M. Levitin, I. Polterovich, and D. A. Sher}{}{\thepage}
}
\usepackage[colorlinks,allcolors=blue,pagebackref]{hyperref}
\usepackage{pifont}
\renewcommand*{\backrefalt}[4]{%
\ifcase #1 %
No citations%
\or
\ding{43}~p.~#2%
\else
\ding{43}~pp.~#2%
\fi}
\newcommand{\mydoi}[1]{\href{https://doi.org/#1}{doi: #1}}

\usepackage{marginnote}

\begin{document}
\pagestyle{mystyle}
\thispagestyle{plain}
\title{%
Inequalities \`a la P\'olya for the Aharonov--Bohm eigenvalues of the disk%
\footnote{A \texttt{Mathematica} script used for a computer-assisted part of the paper and its printout are available for download at \url{https://michaellevitin.net/polya.html\#AB}}
}
\author{
Nikolay Filonov
\thanks{%
\textbf{N. F.: }St. Petersburg Department
of Steklov Institute of Mathematics of RAS,
Fontanka 27, 191023, St. Petersburg, Russia;
St. Petersburg State University,
University emb. 7/9,
199034, St. Petersburg, Russia; 
\href{mailto:filonov@pdmi.ras.ru}{filonov@pdmi.ras.ru}%
}
\and
Michael Levitin\hspace{-3ex}
\thanks{%
\textbf{M. L.: }Department of Mathematics and Statistics, University of Reading, 
Pepper Lane, Whiteknights, Reading RG6 6AX, UK;
\href{mailto:M.Levitin@reading.ac.uk}{M.Levitin@reading.ac.uk}; \url{https://www.michaellevitin.net}%
}
\and 
Iosif Polterovich
\thanks{%
\textbf{I. P.: }D\'e\-par\-te\-ment de math\'ematiques et de statistique, Univer\-sit\'e de Mont\-r\'eal, 
CP 6128 succ Centre-Ville, Mont\-r\'eal QC  H3C 3J7, Canada;
\href{mailto:iossif@dms.umontreal.ca}{iossif@dms.umontreal.ca}; \url{https://www.dms.umontreal.ca/\~iossif}
}
\and
David A. Sher
\thanks{%
\textbf{D. A. S.:  }Department of Mathematical Sciences, DePaul University, 2320 N. Kenmore Ave, 60614, Chicago, IL, USA;
\href{mailto:dsher@depaul.edu}{dsher@depaul.edu}
}
}
\date{{\normalfont \normalsize\em Dedicated to E. Brian Davies on the occasion of his eightieth birthday: to a wonderful man and a brilliant mathematician who taught us, among other things, that proofs and numerics can coexist}\\[1ex]\small arXiv:2311.14072v2, to appear in J. Spect. Theory}
\maketitle

\begin{abstract}  
We prove  an analogue of  P\'olya's conjecture for the eigenvalues of the magnetic Schr\"odinger operator with Aharonov--Bohm potential on the disk, for Dirichlet and magnetic Neumann boundary conditions.
This answers a question posed  by R. L. Frank and A. M. Hansson in 2008.
\end{abstract}

\section{Introduction and main results}\label{sec:intro}

Consider a magnetic Schr\"odinger operator on $\mathbb{R}^2$ with an Aharonov--Bohm potential,
\begin{equation}
\label{eq: AB}
\mathcal{A}_{p,\alpha} := \left(-\mathrm{i}\nabla - \alpha \mathbf{A}_p\right)^2,
\end{equation}
where $p=(p_1,p_2) \in \mathbb{R}^2$, $\alpha \in\mathbb{R}$,  and 
$\mathbf{A}_{p}=|x-p|^{-2}\begin{pmatrix} p_2-x_2\\x_1-p_1\end{pmatrix}$, $x=(x_1, x_2) \in \mathbb{R}^2$.
We refer to $\A_{p,\alpha}$  for brevity as the {\em  Aharonov--Bohm operator} \cite{FH}. It can be rewritten as 
\[
 \A_{p,\alpha} u = -\Delta u + \alpha^2 |\mathbf{A}_p|^2 u + 2\mathrm{i} \alpha \left\langle \nabla u, \mathbf{A}_p\right\rangle,
\]
 where $\Delta=\partial_{x_1}^2 + \partial_{x_2}^2$ is the Laplacian and $u$  is 
 a complex-valued function. The point  $p$ is called the {\em  pole} of the potential,  
and the number $\alpha$ is called the {\em  flux},  cf. \eqref{eq:flux} below. If $\alpha-\alpha' \in \mathbb{Z}$, one can show that the operator $\A_{p,\alpha}$ is unitarily equivalent to $\A_{p, \alpha'}$ (this is also known as {\em  gauge invariance}, see, for instance, \cite[Appendix A.4]{CPS}). Moreover, changing the sign of the flux is equivalent to 
choosing a coordinate system centred at $p$ and interchanging the coordinates. Therefore, we may assume without loss of generality that $\alpha \in \left[0,\frac{1}{2}\right]$.

The operator $\A_{p,\alpha}$  arises in the study of an important phenomenon in quantum physics called the {\em  Aharonov--Bohm effect} \cite{AB}, 
which, roughly speaking, shows that a magnetic potential may affect a charged particle even if  the corresponding magnetic field vanishes. 
Let  
\[
A_p=\frac{p_2-x_2}{|x-p|^2}  \dr x_1 + \frac{x_1-p_1}{|x-p|^2} \dr x_2 
\]
be the dual $1$-form to the vector potential $\mathbf{A}_p$.   Recall that a $1$-form $A$ is called {\em  closed} if $\dr A=0$,  and it is called {\em  exact}
if there is a smooth function $f$ such that $A=\dr f$; it is well-known that adding a magnetic potential corresponding to an exact $1$-form does not change the magnetic Laplacian up to unitary equivalence. One can  check that $A_p$ is a closed form on $\mathbb{R}^2\setminus\{p\}$, 
and hence the magnetic field  $\dr A_p$ generated by the Aharonov--Bohm potential on $\mathbb{R}^2\setminus\{p\}$ is equal to zero.  At the same time, a direct calculation yields that
$A_p = \dr\theta$, where $\theta$ is the angular coordinate in a polar coordinate system centred at $p$.  Therefore, for any simple closed curve $\Gamma$ containing $p$ in its interior,
\begin{equation}
\label{eq:flux}
\frac{1}{2\pi} \oint_\Gamma  A_p = 1,
\end{equation}
and thus the form $\mathcal{A}_p$ is {\em  not} exact. Note that  \eqref{eq:flux} implies that the (normalised) flux of the potential $\alpha \mathbf{A}_p$ is equal to $\alpha$, which motivates the terminology. As shown below, the spectral properties of the operator $\A_{p,\alpha}$ depend on the position of the pole $p$ and  on the value of $\alpha$, despite the fact that the magnetic field is zero.  This can be viewed as a manifestation of the Aharonov--Bohm effect.   There is  a vast literature in spectral theory on the subject, see, for example,  \cite{Hel88, FH, Len,  Bon, CPS} and references therein.

More specifically, let $\Omega \subset \mathbb{R}^2$ be a bounded planar domain.
Consider the eigenvalue problem $\A_{p,\alpha} u = \lambda u$ in $\Omega$ for the Aharonov--Bohm operator with either {\em Dirichlet}
\[
u|_{\partial \Omega}=0
\]
or {\em magnetic Neumann} 
\[
 \left\langle\left.\left(\nabla u-\ir \alpha \mathbf{A}_p  u\right)\right|_{\partial \Omega}, \mathbf{n}\right\rangle = 0
\]
 boundary conditions. Here $\mathbf{n}$ denotes the external unit normal at  the boundary.
 It is well-known that the Dirichlet problem has discrete spectrum for any bounded domain $\Omega$, and the same is true for the 
 magnetic Neumann problem under the standard  regularity assumptions, e.g., if $\partial \Omega$ is Lipschitz, see \cite[Appendix A]{CPS}.
 
Applying a coordinate change if necessary, we may assume without loss of generality that  the pole $p$ coincides with  the origin $o$, and we set 
$\A_\alpha:=\A_{o,\alpha}$, $\alpha \in \left[0,\frac{1}{2}\right]$. 
Moreover, we may assume that $o$  lies in  the simply-connected hull of $\Omega$, i.e., the smallest simply-connected set containing $\Omega$. Indeed,  otherwise the flux of the potential over any closed curve  contained in $\Omega$ is equal to zero, hence the form $A_p$ is exact in $\Omega$,  and the magnetic potential vanishes after an application of a gauge transformation.
 
 We will denote the Dirichlet eigenvalues of $\A_{\alpha}$  on $\Omega$ by 
\[
\lambda^\Dir_1(\Omega,\alpha)\le \lambda^\Dir_2(\Omega,\alpha)\le \dots \nearrow +\infty,
\]
and the magnetic Neumann eigenvalues by 
\[
\lambda^\Neu_1(\Omega,\alpha)\le \lambda^\Neu_2(\Omega,\alpha)\le \dots \nearrow +\infty.
\]
For $\alpha=0$ they coincide, respectively,  with the Dirichlet and Neumann eigenvalues of the usual Laplacian on  $\Omega$.
 
Consider the Dirichlet problem first. Let  
\[
\mathcal{N}_{\Omega}^{\Dir}(\lambda;\alpha) = \#\left\{n\in\mathbb{N}:\lambda^\Dir_n(\Omega,\alpha)\le\lambda^2\right\}
\]
be the eigenvalue counting function. It is known to satisfy Weyl's law (\cite[Theorem A.1]{F08}, see also \cite[Theorem 6]{F22})
\begin{equation}
\label{eq:Weyldir}
\mathcal{N}_{\Omega}^{\Dir}(\lambda;\alpha) =\frac{\area(\Omega)\, \lambda^2}{4\pi}+o(\lambda^2),
\end{equation}
as $\lambda \to \infty$.
If $\alpha=0$, the celebrated P\'olya's eigenvalue conjecture \cite{Pol54} states that for any bounded domain $\Omega \subset \mathbb{R}^2$, 
\begin{equation}
\label{eq:Polyaconj}
\mathcal{N}_{\Omega}^{\Dir}(\lambda;0) \le  \frac{\area(\Omega)\, \lambda^2}{4\pi} 
\end{equation}
for all $\lambda \ge 0$. It was  proved by P\'olya  for tiling domains \cite{Pol}, and has been recently obtained for the disk in \cite{FLPS}; we refer to \cite[\S3.3.3]{LMP} for an overview of this topic.
One may ask if P\'olya's conjecture can be extended to magnetic Schr\"odinger operators. It is known to be false in the case of constant magnetic field \cite{FLoW}. At the same time, numerical experiments presented in \cite{FH} suggested that it holds for the Aharonov--Bohm operator on certain domains, in particular, for the disk. The goal of this paper is to confirm the latter, developing the approach proposed in \cite{FLPS}.   

From now on, let  $\Omega$ be the unit disk $\mathbb{D}$ centred at the origin. 
The (non-ordered) eigenvalues of the Aharonov--Bohm operator $\A_\alpha$ on $\mathbb{D}$ are given by \cite[\S5.2]{FH}
\begin{equation}
\label{eq:ABeigdisk}
\lambda_{i,k}^\Dir:=j_{|i-\alpha|,k}^2,\qquad i\in\mathbb{Z}, k\in\mathbb{N},
\end{equation}
where $j_{\nu,k}$ is the $k$th positive root of the Bessel function $J_\nu(x)$. 

Our first main result is
\begin{theorem}
\label{thm:polyad} Let $\alpha\in\left[0,\frac{1}{2}\right]$. Then
\begin{equation}\label{eq:polyadir}
\mathcal{N}_\disk^{\Dir}(\lambda;\alpha)< \frac{\lambda^2}{4}\qquad\text{for all }\lambda\in(0,+\infty).
\end{equation}
\end{theorem}
We recall \cite[Remark 1.1]{FLPS} that \eqref{eq:polyadir} can be equivalently rewritten as 
\[
\lambda^\Dir_n(\disk,\alpha)>4n\qquad\text{for all }n\in\mathbb{N}.
\]
We prove Theorem \ref{thm:polyad} in \S\ref{sec:proofpd}.

\begin{figure}[htpb]
\centering
\includegraphics{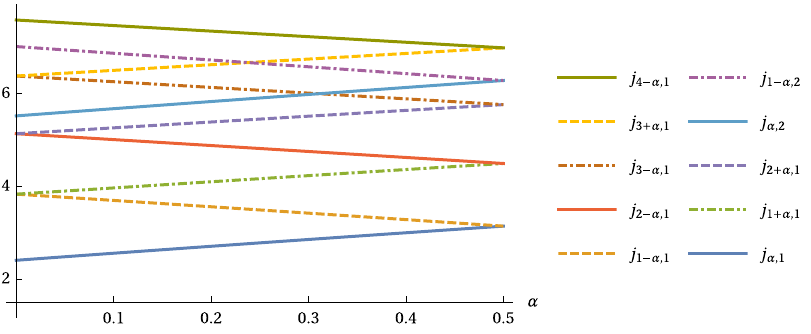}
\caption{Some zeros of Bessel functions as functions of $\alpha$.
\label{fig:bzeros}}
\end{figure}

\begin{remark}
We note that there is no uniform (in $n$) monotonicity relations in $\alpha$ between 
$\lambda_n^\Dir(\disk,\alpha)$ and  $\lambda_n^\Dir(\disk,0)$, and hence \eqref{eq:polyadir} cannot be easily deduced from
the already known case $\alpha=0$.
Indeed, for $\alpha\in\left(0,\frac{1}{2}\right]$ and $m>0$, each double eigenvalue $j_{m,k}^2$ of the Dirichlet Laplacian splits into a smaller eigenvalue $j_{m-\alpha,k}^2$ and a bigger eigenvalue $j_{m+\alpha,k}^2$ of the Aharonov--Bohm operator, whereas the single eigenvalues $j_{0,k}^2$ of the Dirichlet Laplacian  all move up to become the eigenvalues $j_{\alpha,k}^{2}$ of  $\A_\alpha$, see Figure~\ref{fig:bzeros}.
\end{remark}

Consider now the magnetic Neumann problem.  Let 
\[
\mathcal{N}_{\Omega}^{\Neu}(\lambda;\alpha) = \#\left\{n\in\mathbb{N}:\lambda^\Neu_n(\Omega,\alpha)\le\lambda^2\right\}.
\]
For $\alpha=0$ we recover the counting function of the Neumann Laplacian. In this case, 
P\'olya's conjecture states that
\begin{equation}
\label{eq:Polyaconjn}
\mathcal{N}_{\Omega}^{\Neu}(\lambda;0) \ge  \frac{\area(\Omega)\, \lambda^2}{4\pi}.
\end{equation}
Similarly to \eqref{eq:Polyaconj}, inequality \eqref{eq:Polyaconjn} is known to hold for tiling domains \cite{Pol, Kel} and for the disk \cite{FLPS}. We would like to extend it for the Aharonov--Bohm operator on the disk with an arbitrary flux $\alpha$. We note, first of all, that since we assume that the potential has the pole at the origin, the magnetic  Neumann boundary conditions simplify to the ordinary Neumann conditions 
\[
\left\langle \left.\nabla u\right|_{\partial\disk},\mathbf{n}\right\rangle=0.
\]
The (non-ordered) eigenvalues of the Neumann Aharonov--Bohm operator
$\A_\alpha$ on the disk are given
by \cite[\S B.3]{CPS}
\begin{equation}
\label{eq:eigneumann}
\lambda_{i,k}^\Neu:=\left(j'_{|i-\alpha|,k}\right)^2,\qquad i\in\mathbb{Z}, k\in\mathbb{N},
\end{equation}
where $j'_{\nu,k}$ is the $k$th positive root of the derivative of the Bessel function $J'_\nu(x)$, with an exception $j'_{0,1}=0$. 
We immediately note that unlike the case $\alpha=0$ for which \eqref{eq:polyaneu} holds for all $\lambda>0$,  we need to impose some restrictions when $\alpha>0$. Indeed, in this case the first eigenvalue of the Neumann Aharonov--Bohm operator is no longer zero,
\[
\lambda^\Neu_1(\disk,\alpha)=\left(j'_{\alpha,1}\right)^2>0.
\]
Therefore $\mathcal{N}_{\Omega}^{\Neu}(\lambda;0)=0$ for all $\lambda\in\left[0,j'_{\alpha,1}\right)$, and \eqref{eq:Polyaconjn} cannot hold in this interval. However, outside of this interval  P\'olya's conjecture for the Aharonov--Bohm operator on the disk with magnetic Neumann boundary condition still holds.
\begin{theorem}
\label{thm:polyan} Let $\alpha\in\left[0,\frac{1}{2}\right]$. Then 
\begin{equation}\label{eq:polyaneu}
 \mathcal{N}_{\disk}^{\Neu}(\lambda;\alpha)  >\frac{\lambda^2}{4}
\end{equation}
 holds for all $\lambda\ge j'_{\alpha,1}$.
\end{theorem}
Equivalently, we have 
\begin{equation}\label{eq:polyaneun}
\lambda^\Neu_{n+1}(\Omega,\alpha)<4n\qquad\text{for all }n\in\mathbb{N}.
\end{equation}

Theorem \ref{thm:polyan} is proved in \S\ref{sec:Neu}. We note that in the interval $\lambda \in \left[\frac{5}{2},9\right]$
we use a rigorous computer-assisted argument based on integer arithmetic, see Proposition \ref{prop:gap}.

\begin{remark}
In fact, apart from the ground state \cite{CS, CPS},  the magnetic Neumann eigenvalues of the Aharonov--Bohm operator 
appear to be  less explored than the  Dirichlet ones. In particular, even the one-term  Weyl's law,
\begin{equation}
\label{eq:Weylneu}
\mathcal{N}_{\Omega}^{\Neu}(\lambda;\alpha) =\frac{\area(\Omega)\, \lambda^2}{4\pi}+o(\lambda^2)\qquad\text{as }\lambda\to+\infty,
\end{equation}
has not been worked out in the literature in full detail. We note that the standard smooth methods do not apply due the presence of a singularity at the pole. To fill in this gap, we present an argument
that was communicated to us by R. L. Frank \cite{F23}. As before, assume that $\Omega$ has Lipschitz boundary and consider the Aharonov--Bohm operator $\A_\alpha$ with a pole at the origin $o \in \Omega$. 
Since $\mathbf{A}_o \in L^2_{\mathrm{loc}} \left(\Omega \setminus \{o\}\right)$,  it follows from \cite[Theorem 1.3]{HS04} that the pointwise difference between the Neumann and Dirichlet heat kernels for $\A_\alpha$ on the diagonal satisfies
\[
\left| \er^{- t \A_{\alpha}^\Neu}(x,x) - \er^{- t \A_{\alpha}^\Dir} (x,x) \right| \leq \er^{t \Delta^\Neu} (x,x) - \er^{t \Delta^\Dir} (x,x), \quad  x \in \Omega \setminus \{o\}, \quad t>0.
\] 
Integrating both sides over $x$, multiplying by $t^{d/2}$ and using the heat trace asymptotics for the Dirichlet and Neumann Laplacian, we get
\begin{equation}\label{eq:heat}
t^{d/2} \left(\operatorname{Tr} \er^{-t \A_{\alpha}^\Neu} - \operatorname{Tr} \er^{- t \A_{\alpha}^\Dir}\right)  \to 0\qquad\text{as }t\to 0^+.
\end{equation}
By an Abelian theorem, the Weyl asymptotics \eqref{eq:Weyldir} implies the one-term asymptotics for the Dirichlet heat trace for the Aharonov--Bohm operator, which by \eqref{eq:heat} implies the same asymptotics for the Neumann heat trace. The standard Tauberian argument then implies  \eqref{eq:Weylneu}.
%
\end{remark}

\begin{remark}
It would be interesting to know if analogues of Theorems \ref{thm:polyad} and \ref{thm:polyan} hold for other domains, even numerically.  As we have mentioned previously, in the Dirichlet case it was checked computationally in \cite{FH} for some annuli and the square with the centred potential. We have conducted some rough numerical experiments for the disk, with both the Dirichlet and the Neumann conditions, and an off-centred pole position, which seem to indicate that \emph{experimentally} the analogues of P\'olya's conjecture still hold in these cases (with an appropriate modification for the first Neumann eigenvalue as above) for all $\alpha$ and all pole placements, but we do not have a proof of this claim. 
\end{remark}

\subsection*{Acknowledgements}  We are grateful to Rupert Frank and Jens Marklof for drawing our attention to this problem. We would also like to thank Rupert Frank for many useful discussions and for communicating to us the proof of \eqref{eq:Weylneu}.  Research of N. F. was supported by the grant No.\ 22-11-00092 of the Russian Science Foundation.  Research of M. L. was partially supported by the EPSRC and by the University of Reading RETF Open Fund.  Research of I. P. was partially supported by NSERC. The authors would also like to acknowledge the hospitality of the University of Bristol and the Mathematisches Forschungsinstitut Oberwolfach  where this work has been conceived.

\section{Proof of Theorem {\ref{thm:polyad}}} 
\label{sec:proofpd}
Recall formula \eqref{eq:ABeigdisk}.
With an alternative parametrisation, we can split the eigenvalues into two doubly-indexed sequences,
\[
\begin{aligned}
&j_{m+\alpha, k}^2,\\
&j_{m+(1-\alpha), k}^2,
\end{aligned}
\qquad m\in\mathbb{N}\cup\{0\}, k\in\mathbb{N}.
\]
Therefore, the counting function of eigenvalues which are less than or equal to a given $\lambda^2\in[0,+\infty)$ is 
\begin{equation}\label{eq:Nlambda}
\begin{split}
\mathcal{N}^\Dir_\disk(\lambda;\alpha)&:=\#\left\{n\in\mathbb{N}:\lambda^\Dir_n(\disk,\alpha)\le\lambda^2\right\}=\#\left\{(i,k)\in \mathbb{Z}\times \mathbb{N}: \lambda_{i,k}^\Dir\le \lambda^2\right\}\\
&=\sum_{m=0}^\infty \#\left\{k\in \mathbb{N}: j_{m+\alpha, k}\le\lambda\right\}+\sum_{m=0}^\infty \#\left\{k\in \mathbb{N}: j_{m+(1-\alpha), k}\le\lambda\right\}\\
&=\sum_{m=0}^{\entire{\lambda-\alpha}} \#\left\{k\in \mathbb{N}: j_{m+\alpha, k}\le\lambda\right\}+\sum_{m=0}^{\entire{\lambda-(1-\alpha)}} \#\left\{k\in \mathbb{N}: j_{m+(1-\alpha), k}\le\lambda\right\},
\end{split}
\end{equation}
where in the last line we used the fact that $j_{\nu,1}> \nu$, and therefore $\#\left\{k\in \mathbb{N}: j_{\nu, k}\le\lambda\right\}=0$ whenever $\nu>\lambda$. Here and further on, $\entire{x}$ and $\ceiling{x}$ stand, as usual, for the floor and the ceiling of a real number $x$, respectively.

The main ideas of the proof are similar to those in the case of the Dirichlet Laplacian \cite{FLPS}, with some notable changes. The first important ingredient is a uniform upper bound on the number of zeros of a Bessel function in a given interval, which is directly borrowed from  \cite{FLPS} and \cite{Sh}. For $\lambda\ge 0$, we define the function $G_\lambda:[0,+\infty)\to \left[0,\frac{\lambda}{\pi}\right]$ as
\begin{equation}\label{eq:G}
G_\lambda(z):=\begin{cases}
\frac{1}{\pi}\left(\sqrt{\lambda^2-z^2}-z\arccos\frac{z}{\lambda}\right)&\qquad\text{if }z\in[0,\lambda],\\
0&\qquad\text{if }z\in(\lambda,+\infty).
\end{cases}
\end{equation}
We note that $G_\lambda(z)$ is convex on its domain, strictly monotone decreasing on $[0,\lambda]$ with $-\frac{1}{2}\le G'_\lambda(z)\le 0$ everywhere, and that 
\begin{equation}\label{eq:intG}
\int_0^\lambda G_\lambda(z)\,\dr z=\frac{\lambda^2}{8},
\end{equation}
see \cite[Lemmas 4.5 and 4.6]{FLPS}.

\begin{lemma}[{\cite[Proposition 3.1]{FLPS}}]\label{lem:bes} Let $\nu\ge 0$ and $\lambda> 0$. Then
\[
\#\left\{k: j_{\nu,k}\le \lambda\right\}\le \entire{G_\lambda(\nu)+\frac{1}{4}}.
\]
\end{lemma}
 
 Applying now Lemma \ref{lem:bes} to the right-hand side of equation \eqref{eq:Nlambda}, we immediately get
 \begin{equation}\label{eq:Palpha}
\mathcal{N}^\Dir_\disk(\lambda;\alpha) \le \sum_{m=0}^{\ceiling{\lambda}-1} \left( \entire{G_\lambda(m+\alpha)+\frac{1}{4}}+\entire{G_\lambda\left(m+(1-\alpha)\right)+\frac{1}{4}}\right)=:\mathcal{P}^\Dir_\alpha(\lambda).
 \end{equation}
 
\begin{remark} It may happen that for $m=\ceiling{\lambda}-1$, the argument $m+(1-\alpha)$, either on its own or together with $m+\alpha$, exceeds $\lambda$, in which case the corresponding terms do not contribute to $\mathcal{P}^\Dir_\alpha(\lambda)$ according to the definition of $G_\lambda$.
\end{remark} 

 We note that the quantity $\mathcal{P}^\Dir_\alpha(\lambda)$ appearing in the right-hand side of \eqref{eq:Palpha} is exactly the number of those points of shifted lattices $\left(\mathbb{Z}\pm\alpha\right)\times\left(\mathbb{Z}-\frac{1}{4}\right)$ which lie in the first quadrant under or on the graph of the function $G_\lambda$.  
 
 Some novel methods for estimating the number of shifted lattice points under a graph of a convex decreasing function have been developed in  \cite{FLPS}; they are not, however, directly applicable to the sum $\mathcal{P}^\Dir_\alpha(\lambda)$ for $\alpha>0$. The case $\alpha=0$ has been dealt with in \cite[Theorem 5.1]{FLPS}; the following result is its extension for $\alpha\in\left(0, \frac{1}{2}\right]$. This is the second main ingredient of the proof of Theorem \ref{thm:polyad}.

\begin{theorem}\label{thm:Pcount} Let $b\in\mathbb{N}$ and $\alpha\in\left[0,\frac{1}{2}\right]$. Let $g$ be a non-negative decreasing convex function on $[0, b]$ such that $g(b)=0$ and
\begin{equation}\label{eq:Lip}
|g(z)-g(w)|\le \frac{1}{2}|z-w|
\end{equation}
for all $z,w\in[0,b]$.
Then 
\begin{equation}\label{eq:gsum}
\sum_{m=0}^{b-1} \left( \entire{g\left(m+\alpha\right)+\frac{1}{4}}+\entire{g\left(m+(1-\alpha)\right)+\frac{1}{4}}\right)\le 2\int_0^b g(z)\,\dr z.
\end{equation}
\end{theorem}

The proof of  Theorem \ref{thm:Pcount} is based on the following 

\begin{lemma}\label{lem:g}
Let $i, j\in\mathbb{Z}$, $i<j$, let $\alpha\in\left[0,\frac{1}{2}\right]$, and let $g$ be a decreasing convex function on $[i,j+1]$ satisfying \eqref{eq:Lip} for all $z,w\in  [i,j+1]$. We assume that
\begin{equation}\label{eq:gn}
n+1\ge g(i+1)\ge g(j) \ge n \ge g(j+1),
\end{equation}
and that 
\begin{equation}\label{eq:gnX}
Z:=\min\left\{z\in[i,j+1]: g(z)=n\right\}<j+1.
\end{equation}
Then 
\begin{equation}\label{eq:gsumn}
\sum_{m=i}^{j-1} \left( \entire{g\left(m+\alpha\right)+\frac{1}{4}}+\entire{g\left(m+(1-\alpha)\right)+\frac{1}{4}}\right)\le 2\int_i^j g(z)\,\dr z.
\end{equation}
\end{lemma}

\begin{proof}[Proof of Lemma {\ref{lem:g}}] The validity of the statement does not change if an integer constant is added to $g$, so without loss of generality we can take $n=0$, in which case  \eqref{eq:gn} becomes
\[
1\ge g(i+1)\ge g(j) \ge 0 \ge g(j+1),
\]
and  \eqref{eq:gnX} becomes
\[
Z:=\min\{z\in[i,j+1]: g(z)=0\}<j+1,
\]
which means that we exclude the situation when $x=j+1$ is the first and only zero of $g(x)$ in the interval $[i,j+1]$. 
 
Consider the sequence 
\[
g(i+\alpha)\ge g(i+(1-\alpha))\ge \dots \ge g(j-1+\alpha)\ge  g(j-\alpha). 
\]
We will distinguish five cases depending on which element of this sequence first drops below $\frac{3}{4}$.

\begin{description}
\item[Case 1:  $\frac{3}{4}>g(i+\alpha)$.] The left-hand side of \eqref{eq:gsumn} is zero, and the right-hand side is nonnegative. 
\item[Case 2:  $g(i+\alpha)\ge \frac{3}{4}>g(i+1-\alpha)$.]  The left-hand side of \eqref{eq:gsumn} equals one. We have $g(i)\ge g(i+\alpha)\ge \frac{3}{4}$ by monotonicity, and therefore $g(z)\ge \frac{3}{4}-\frac{z-i}{2}$ for $z\in[i,j]\supseteq[i,i+1]$ by  \eqref{eq:Lip}. Thus,
\[
2\int_i^j g(z)\,\dr z\ge 2\int_i^{i+1} \left(\frac{3}{4}-\frac{z-i}{2}\right)\,\dr z=1.
\]
\item[Case 3:  $g(i+1-\alpha)\ge \frac{3}{4}>g(i+1+\alpha)$.] The left-hand side of \eqref{eq:gsumn} equals two. We have $g\left(i+\frac{1}{2}\right)\ge g(i+1-\alpha)\ge \frac{3}{4}$ by monotonicity, and therefore $g(i+1)\ge \frac{1}{2}$ and  $g(i+2)\ge 0$ by \eqref{eq:Lip}. This is also ensures that $Z\ge i+2$, hence $[i,j]\supseteq[i,i+2]$, and by convexity 
\[
2\int_i^j g(z)\,\dr z\ge 2\int_i^{i+2}g(z)\,\dr z\ge 4g(i+1)\ge 2.
\]
\item[Case 4:  $g(i+k+\alpha)\ge \frac{3}{4}>g(i+k+1-\alpha)$ for some $k\in\mathbb{N}$.] The left-hand side of \eqref{eq:gsumn} is equal to $2k+1$. We have, by convexity,
\[
g(i+1)+g(i+2k+2\alpha-1)\ge 2g(i+k+\alpha)\ge \frac{3}{2},
\]
and therefore $g(i+2k+2\alpha-1)\ge \frac{3}{2}-g(i+1)\ge \frac{1}{2}$.  Then  by  \eqref{eq:Lip}, $g(i+2k)\ge \frac{1-(1-2\alpha)}{2}\ge 0$, and additionally $Z\ge i+2k$, so that $[i,j]\supseteq[i,i+2k]$, and by convexity and monotonicity,
\[
2\int_i^j g(z)\,\dr z\ge 2\int_i^{i+2k}g(z)\,\dr z\ge 4kg(i+k)\ge 4kg(i+k+\alpha)\ge 3k\ge 2k+1
\]
as $k\ge 1$.
\item[Case 5:  $g(i+k-\alpha)\ge \frac{3}{4}>g(i+k+\alpha)$ for some integer $k\ge 2$.] The left-hand side of \eqref{eq:gsumn} is equal to $2k$. By convexity,
\[
g(i+1)+g(i+2k-2\alpha-1)\ge 2g(i+k-\alpha)\ge \frac{3}{2},
\]
and therefore $g(i+2k-2\alpha-1)\ge \frac{3}{2}-g(i+1)\ge \frac{1}{2}$. Then  by  \eqref{eq:Lip}, $g(i+2k-1)\ge \frac{1-2\alpha}{2}\ge 0$, and $Z\ge i+2k-1$, so that $[i,j]\supseteq[i,i+2k-1]$, and by convexity and monotonicity,
\[
\begin{split}
2\int_i^j g(z)\,\dr z\ge 2\int_i^{i+2k-1}g(z)\,\dr z&\ge 2(2k-1)g\left(i+k-\frac{1}{2}\right)\\&\ge 2(2k-1)g(i+k-\alpha)\ge \frac{3(2k-1)}{2}> 2k
\end{split}
\]
as $k\ge 2$.
\end{description}
\end{proof}

\begin{remark}\label{rem:linear}
It is easily checked that the equality in \eqref{eq:gsumn} can be attained only if the function $g$ is linear.
\end{remark}

We can now proceed to the
\begin{proof}[Proof of Theorem \ref{thm:Pcount}] We extend the function $g$ by zero onto the interval $(b,b+1]$. Let $N=\entire{g(0)}\ge 0$. We have two cases.
\begin{description}
\item[Case 1: $N=0$.] In this case 
\[
1> g(0)\ge\dots \ge g(b)=0.
\]
We apply Lemma \ref{lem:g} with $i=0$ and $j=b$ (and hence $j+1>b\ge Z$), immediately yielding \eqref{eq:gsum}.
\item[Case 2: $N\ge 1$.]  We set, for $k=0,1,\dots,N$, 
\[
L_k:=\max\{m\in\{0,\dots,b\}: g(m)\ge k\},
\]
see \cite[Figure 5]{FLPS} for an illustration. Since, by the conditions of the theorem,  the function $g$ is strictly monotone decreasing at all points where it takes positive values, the inverse function $g^{-1}$ is well-defined on $(0,g(0)]$, and $L_k=\entire{g^{-1}(k)}$ for $k>0$; also, by definition, $L_0=b$. 
We also set
\[
N_1=\begin{cases}
N&\qquad\text{if }L_N=0,\\
N+1&\qquad\text{if }L_N>0,
\end{cases}
\]
and in the latter case set $L_{N+1}:=0$. 
In any case, we have $0=L_{N_1}<\dots<L_0=b$, where the strict inequalities follow from \eqref{eq:Lip}.

Let $n\in\{N_1-1,\dots,0\}$, and consider the interval $[L_{n+1}, L_n]$. If $n>0$, then, by construction, we have 
\[
n+1>g\left(L_{n+1}+1\right)>\dots>g\left(L_{n}\right)\ge n> g\left(L_{n}+1\right).
\]
Similarly, if $n=0$, then
\[
1>g\left(L_{1}+1\right)>\dots>g\left(L_0\right)= 0=g(b+1).
\]
In either case, we can apply Lemma \ref{lem:g} with $i=L_{n+1}$ and $j=L_n$, giving
\begin{equation}\label{eq:sumgint}
\sum_{m=L_{n+1}}^{L_n-1} \left( \entire{g\left(m+\alpha\right)+\frac{1}{4}}+\entire{g\left(m+(1-\alpha)\right)+\frac{1}{4}}\right)\le 2\int_{L_{n+1}}^{L_n} g(z)\,\dr z.
\end{equation}
Summing up \eqref{eq:sumgint} over $n\in\{N_1-1,\dots,0\}$ with account of $L_{N_1}=0$ and $L_0=b$ yields \eqref{eq:gsum}.
\end{description}
\end{proof}

The result of Theorem \ref{thm:polyad} now follows immediately from applying Theorem \ref{thm:Pcount} to $g=G_\lambda$ with $b=\ceiling{\lambda}$, using \eqref{eq:intG} and \eqref{eq:Palpha}, and recalling that $G_\lambda(z)\equiv 0$ for $z\in\left[\lambda, \ceiling{\lambda}\right]$. The strict inequality follows from Remark \ref{rem:linear} and the fact that $G_\lambda$ is not (piecewise) linear.

\section{Proof of Theorem \ref{thm:polyan}}\label{sec:Neu}

\begin{figure}[htpb]
\centering
\includegraphics{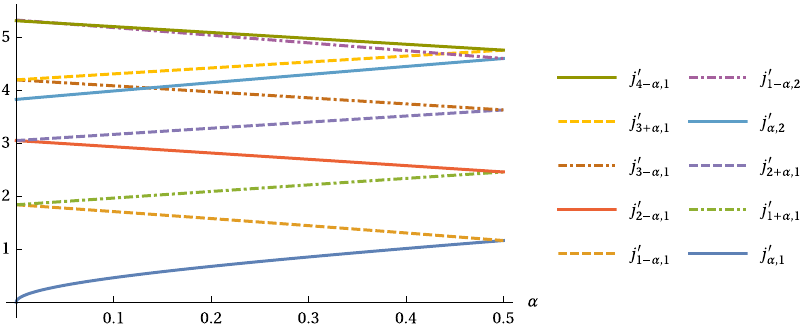}
\caption{Some zeros of derivatives of Bessel functions as functions of $\alpha$.
\label{fig:dbzeros}}
\end{figure}

As in \S\ref{sec:proofpd}, in view of \eqref{eq:eigneumann} the counting function of Neumann eigenvalues which are less than or equal to $\lambda^2\in[0,+\infty)$ is given by
\begin{equation}\label{eq:NlambdaN}
\begin{split}
\mathcal{N}^\Neu_\disk(\lambda;\alpha)&:=\#\left\{n\in\mathbb{N}:\lambda^\Neu_n(\disk,\alpha)\le\lambda^2\right\}=\#\left\{(i,k)\in \mathbb{Z}\times \mathbb{N}: \lambda^\Neu_{i,k}\le \lambda^2\right\}\\
&=\sum_{m=0}^\infty \#\left\{k\in \mathbb{N}: j'_{m+\alpha, k}\le\lambda\right\}+\sum_{m=0}^\infty \#\left\{k\in \mathbb{N}: j'_{m+(1-\alpha), k}\le\lambda\right\}\\
&=\sum_{m=0}^{\entire{\lambda-\alpha}} \#\left\{k\in \mathbb{N}: j'_{m+\alpha, k}\le\lambda\right\}+\sum_{m=0}^{\entire{\lambda-(1-\alpha)}} \#\left\{k\in \mathbb{N}: j'_{m+(1-\alpha), k}\le\lambda\right\},
\end{split}
\end{equation}
where we used the fact that $j'_{\nu,1}\ge \nu$, and therefore $\#\left\{k\in \mathbb{N}: j'_{\nu, k}\le\lambda\right\}=0$ whenever $\nu>\lambda$. The graphs of the first few zeros of Bessel derivatives appearing in \eqref{eq:NlambdaN} are shown in Figure \ref{fig:dbzeros}.

We split the proof of Theorem \ref{thm:polyan} into three parts, for different intervals of $\lambda$. The proofs for very small and sufficiently large $\lambda$s are analytic, and a small remaining gap is closed with the help of a rigorous computer-assisted algorithm similar to that in \cite[\S 8]{FLPS}.

We start with the case of small $\lambda$s. 

\begin{prop}\label{prop:smallambda}
For any $\alpha\in\left(0,\frac{1}{2}\right]$, the bound \eqref{eq:polyaneu} holds for all $\lambda\in\left[j'_{\alpha,1},2\sqrt{2}\right)$.
\end{prop}

\begin{proof}[Proof of Proposition \ref{prop:smallambda}] It is sufficient to show that for all $\alpha\in\left(0,\frac{1}{2}\right]$, \eqref{eq:polyaneun} holds with $n=1$ and $n=2$, namely, that we have
\begin{equation}\label{eq:Nlambda2}
\lambda^\Neu_2(\disk,\alpha)<4
\end{equation}
and
\begin{equation}\label{eq:Nlambda3}
\lambda^\Neu_3(\disk,\alpha)<8.
\end{equation}
It follows immediately from the monotonicity with respect to $\nu$ of the zero $j'_{\nu,k}$  \cite[\S10.21(iv)]{DLMF}, and from the interlacing properties of these zeros \cite[formula (10)]{PA} that
\[
\lambda^\Neu_2(\disk,\alpha)=\left(j'_{1-\alpha,1}\right)^2\qquad\text{and}\qquad\lambda^\Neu_3(\disk,\alpha)=\left(j'_{1+\alpha,1}\right)^2,
\]
cf.\ Figure \ref{fig:dbzeros}. 

As $j'_{1-\alpha,1}\le j'_{1,1}$ and $j'_{1+\alpha,1}\le j'_{\frac{3}{2},1}$ for all $\alpha\in\left[0,\frac{1}{2}\right]$, we, in principle, could finish the proof of Proposition  \ref{prop:smallambda} here by an elementary observation that $j'_{1,1}\approx 1.8412$ and $j'_{\frac{3}{2},1}\approx 2.4605$. However we deliberately avoid the use of floating point approximations of zeros of Bessel functions and their derivatives anywhere in this paper, so we continue to give a fully analytic justification of \eqref{eq:Nlambda2} and \eqref{eq:Nlambda3}.

Since $j'_{1,1}=\sqrt{\lambda^\Neu_2(\disk,0)}$, which is less than $2$ by P\'olya's conjecture for the Neumann Laplacian in the disk, we have \eqref{eq:Nlambda2}.

The Bessel functions of half-integer order can be expressed explicitly in terms of trigonometric functions and powers of their argument,  and as a result $j'_{\frac{3}{2},1}$ coincides with the first positive root of the function 
\[
S(z):=-\cot z -\frac{2}{3}z+\frac{1}{z}.
\]
Elementary calculations show that 
\[
\begin{gathered}
\lim_{z\to 0^+}S(z)=0,\quad \lim_{z\to 0^+}S'(z)=-\frac{1}{3},\quad \lim_{z\to\pi^-}S'(z)=+\infty,\\
S''(z)=2\left(\frac{1}{z^3}-\frac{\cos z}{\sin^3 z}\right)>0\qquad\text{for all }z\in(0,\pi),
\end{gathered}
\]
and therefore $S(z)$ has a single root in $(0,\pi)$. As $S\left(\frac{5\pi}{6}\right)=\sqrt{3}+\frac{6}{5\pi}-\frac{5\pi}{9}>0$, we conclude that 
\[
j'_{\frac{3}{2},1}<\frac{5\pi}{6}<2\sqrt{2},
\]
thus proving \eqref{eq:Nlambda3}.
\end{proof}

We now switch to large values of $\lambda$s and prove
\begin{prop}\label{prop:largelambda}
For any $\alpha\in\left(0,\frac{1}{2}\right]$, the bound \eqref{eq:polyaneu} holds for all $\lambda\ge \frac{2500\pi}{3481\pi-10000}$.
\end{prop}

We note that $\frac{2500\pi}{3481\pi-10000}\approx 8.39205$.

The proof of Proposition \ref{prop:largelambda} is based on several results from \cite{FLPS}, the first of which complements Lemma \ref{lem:bes}.
\begin{lemma}[{\cite[Proposition 3.1]{FLPS}}]\label{lem:dbes} Let $\nu\ge 0$ and $\lambda> 0$. Then
\[
\#\left\{k: j'_{\nu,k}\le \lambda\right\}\ge \entire{G_\lambda(\nu)+\frac{3}{4}}.
\]
\end{lemma}

Applying Lemma \ref{lem:dbes} to the right-hand side of \eqref{eq:NlambdaN}, we get
\begin{equation}\label{eq:PalphaN}
\mathcal{N}^\Neu_\disk(\lambda;\alpha) \ge
\sum_{m=0}^{\entire{\lambda-\alpha}} \entire{G_\lambda(m+\alpha)+\frac{3}{4}}+ \sum_{m=0}^{\entire{\lambda-1+\alpha}}\entire{G_\lambda\left(m+(1-\alpha)\right)+\frac{3}{4}}
=:\mathcal{P}^\Neu_\alpha(\lambda).
\end{equation}

In order to obtain a lower bound on $\mathcal{P}^\Neu_\alpha(\lambda)$, we recall
\begin{lemma}[{\cite[Theorem 6.1]{FLPS}}]\label{lem:sumNeu}
Let $b>0$, and let $g$ be a non-negative decreasing convex function on $[0,b]$ such that $g(0)\ge \frac{1}{4}$, $g(b)=0$, and \eqref{eq:Lip} holds for all $z,w\in[0,b]$. Let
\[
M_0=M_{g,0}:=1+\max\left\{m\in\{0,\dots,\entire{b}\}: g(m)\ge \frac{1}{4}\right\},
\]
and assume that $M_0\le b$. 
Then 
\begin{equation}\label{eq:60}
\sum_{m=0}^{\entire{b}} \entire{g(m)+\frac{3}{4}}\ge \int_0^{M_0} g(z)\,\dr z+\frac{M_0}{4}.
\end{equation}
\end{lemma}

\begin{remark}\label{rem:M0} If $g$ is strictly monotone on $[0, b]$, then
\[
M_0=\entire{g^{-1}\left(\frac{1}{4}\right)}+1,
\]
\end{remark}

\begin{remark} 
The bound in \cite[Theorem 6.1]{FLPS} was stated slightly differently compared to \eqref{eq:60}, giving a marginally weaker result,  however inequality \eqref{eq:60} also appeared there towards  the end of the proof.
\end{remark}

\begin{proof}[Proof of Proposition \ref{prop:largelambda}]
We want to apply Lemma \ref{lem:sumNeu} to each of the sums in the right-hand side of  \eqref{eq:PalphaN}. 

For the first sum, we take $g(z)=g_1(z):=G_\lambda(z+\alpha)$, $b=\lambda-\alpha$, and 
\begin{equation}\label{eq:Mg1}
M_0=M_{g_1,0}:=\entire{g_1^{-1}\left(\frac{1}{4}\right)}+1=\entire{G_\lambda^{-1}\left(\frac{1}{4}\right)-\alpha}+1,
\end{equation} 
and require that $g_1(0)=G_\lambda(\alpha)\ge \frac{1}{4}$. We note also that in this case $M_0=\entire{G_\lambda^{-1}\left(\frac{1}{4}\right)-\alpha}+1\le G_\lambda^{-1}\left(\frac{1}{4}\right)-\alpha+1 \le b=\lambda-\alpha$ if $\lambda\ge 2$ by \cite[Lemma 4.8]{FLPS}, so from now on we assume this restriction on $\lambda$. Applying Lemma \ref{lem:sumNeu} gives
\begin{equation}\label{eq:Sigma1}
\begin{split}
\Sigma_1(\lambda):=\sum_{m=0}^{\entire{\lambda-\alpha}}  \entire{G_\lambda(m+\alpha)+\frac{3}{4}}&\ge \int_0^{M_{g_1,0}} G_\lambda(z+\alpha)\,\dr z+\frac{M_{g_1,0}}{4}\\
&=\int_{\alpha}^{M_{g_1,0}+\alpha} G_\lambda(z)\,\dr z +\frac{M_{g_1,0}}{4}\\
&\ge \frac{\lambda^2}{8}-\left(\int_0^\alpha+\int_{M_{g_1,0}+\alpha}^\lambda\right) G_\lambda(z)\,\dr z +\frac{M_{g_1,0}}{4},
\end{split}
\end{equation}
where we used \eqref{eq:intG}.

For the second sum, we take $g(z):=g_2(z):=G_\lambda(z+1-\alpha)$, $b=\lambda-1+\alpha$, and 
$M_0=M_{g_2,0}:=\entire{g^{-1}\left(\frac{1}{4}\right)}+1=\entire{G_\lambda^{-1}\left(\frac{1}{4}\right)-1+\alpha}+1$, and require that $g_{2}(0)=G_\lambda(1-\alpha)\ge \frac{1}{4}$. Similarly to the above, in this case $M_0\le b$ if $\lambda\ge 2$, and applying  Lemma \ref{lem:sumNeu} we get
\begin{equation}\label{eq:Sigma2}
\begin{split}
\Sigma_2(\lambda):=\sum_{m=0}^{\entire{\lambda-1+\alpha}}  \entire{G_\lambda(m+1-\alpha)+\frac{3}{4}}&\ge \int_0^{M_{g_2,0}} G_\lambda(z+1-\alpha)\,\dr z+\frac{M_{g_2,0}}{4}\\
&\ge \frac{\lambda^2}{8}-\left(\int_0^{1-\alpha}+\int_{M_{g_2,0}+1-\alpha}^\lambda\right) G_\lambda(z)\,\dr z +\frac{M_{g_2,0}}{4},
\end{split}
\end{equation}

Adding the two bounds \eqref{eq:Sigma1} and \eqref{eq:Sigma2} yields
\begin{equation}\label{eq:PNbound}
\begin{split}
\mathcal{P}^\Neu_\alpha(\lambda)&=\Sigma_1(\lambda)+\Sigma_2(\lambda)\\
&\ge  \frac{\lambda^2}{4}-\left(\int_0^{\alpha}+\int_0^{1-\alpha}+\int_{M_{g_1,0}+\alpha}^\lambda+\int_{M_{g_2,0}+1-\alpha}^\lambda\right) G_\lambda(z)\,\dr z+\frac{M_{g_1,0}+M_{g_2,0}}{4}.
\end{split}
\end{equation}

We note that we additionally require
\begin{equation}\label{eq:Gcond}
G_\lambda(\alpha)\ge G_\lambda(1-\alpha)\ge \frac{1}{4}
\end{equation}
for \eqref{eq:Sigma1}--\eqref{eq:PNbound} to be valid.
We now check if we need to restrict the range of $\lambda$ further to ensure that  \eqref{eq:Gcond} holds. We can replace  \eqref{eq:Gcond} by a stronger condition 
\[
G_\lambda(1)=\frac{1}{\pi}\left(\sqrt{\lambda^2-1}-\arccos\frac{1}{\lambda}\right)\ge \frac{1}{4}.
\]
Since $G_\lambda(1)$ is monotone increasing in $\lambda$,  and $G_3(1)>\frac{2\sqrt{2}}{\pi}-\frac{1}{2}>\frac{1}{4}$, it is sufficient to assume from now on that  $\lambda\ge 3$.

To make the bound in the right-hand side of \eqref{eq:PNbound} more explicit in terms of $\lambda$, we first estimate the integrals there. For integrals near zero, an obvious bound 
\begin{equation}\label{eq:intnear0}
\int_0^\beta G_\lambda(z)\,\dr z\le \beta G_\lambda(0)=\frac{\beta\lambda}{\pi}\qquad\text{for all }\beta\in[0,\lambda].
\end{equation}
holds due to monotonicity of $G_\lambda$. To estimate the integrals near $\lambda$, we use the following
\begin{lemma}\label{lem:intGbound} 
Let $G_\lambda$ be defined by \eqref{eq:G}, and let $\beta\in[0,\lambda]$. Then
\[
\int_{\beta}^\lambda G_\lambda(z)\,\dr z\le \frac{2}{5}(\lambda-\beta)G_\lambda(\beta).
\]
\end{lemma}

\begin{proof}[Proof of Lemma \ref{lem:intGbound}]
We fix $\lambda>0$, change the variable $\beta=\lambda\gamma$, $\gamma\in[0,1]$, and intend to show that 
\[
\rho(\gamma):=\int_{\gamma\lambda}^\lambda G_\lambda(z)\,\dr z-\frac{2}{5}\lambda(1-\gamma)G_\lambda(\gamma\lambda)\le 0\qquad\text{for all  }\gamma\in[0,1].
\]
We have, by an explicit calculation, 
\[
\rho'(\gamma)=\frac{(\gamma+2)\lambda^2}{5\pi}\left(\arccos\gamma - \frac{3\sqrt{1-\gamma^2}}{\gamma+2}\right),
\]
with $\rho(1)=0$, therefore it is sufficient to show that 
\[
\chi(\gamma):=\arccos\gamma - \frac{3\sqrt{1-\gamma^2}}{\gamma+2}\ge 0\qquad\text{for all  }\gamma\in[0,1].
\]
We again have $\chi(1)=0$, and the result follows from the fact that
\[
\chi'(\gamma)=-\frac{(1-\gamma)^2}{(\gamma+2)^2 \sqrt{1-\gamma ^2}}\le 0\qquad\text{for all  }\gamma\in[0,1].
\]
\end{proof}

We now return to the previously obtained inequality \eqref{eq:PNbound} and estimate the integrals in its right-hand side using \eqref{eq:intnear0}. This gives
\begin{equation}\label{eq:PNeubound2}
\begin{split}
\mathcal{P}^\Neu_\alpha(\lambda)&\ge  \frac{\lambda^2}{4}-\left(\frac{\lambda}{\pi}+\frac{2}{5}\left(\left(\lambda-M_{g_1,0}-\alpha\right)G_\lambda\left(M_{g_1,0}+\alpha\right)+\left(\lambda-M_{g_2,0}-1+\alpha\right)G_\lambda\left(M_{g_2,0}+1-\alpha\right)\right)\right)\\
&\qquad+\frac{M_{g_1,0}+M_{g_2,0}}{4}\\
&=\frac{\lambda^2}{4}-\lambda\left(\frac{1}{\pi}+\frac{2}{5}\left(G_\lambda\left(M_{g_1,0}+\alpha\right)+G_\lambda\left(M_{g_2,0}+1-\alpha\right)\right)\right)\\
&+\frac{2}{5}\left(\left(M_{g_1,0}+\alpha\right) G_\lambda\left(M_{g_1,0}+\alpha\right)+\left(M_{g_2,0}+1-\alpha\right) G_\lambda\left(M_{g_2,0}+1-\alpha\right)\right)+\frac{M_{g_1,0}+M_{g_2,0}}{4}.
\end{split}
\end{equation}
From the definition \eqref{eq:Mg1} of $M_{g_1,0}$, the inequality $\entire{x}+1\ge x$, and monotonicity of $G_\lambda$ we obtain
\[
M_{g_1,0}+\alpha\ge G_\lambda^{-1}\left(\frac{1}{4}\right),\qquad G_\lambda\left(M_{g_1,0}+\alpha\right)\le \frac{1}{4},
\]
and exactly in the same manner 
\[
M_{g_2,0}+1-\alpha\ge G_\lambda^{-1}\left(\frac{1}{4}\right),\qquad G_\lambda\left(M_{g_2,0}+1-\alpha\right)\le  \frac{1}{4}.
\]
Substituting these bounds into \eqref{eq:PNeubound2}, we get, after simplifications,
\begin{equation}\label{eq:PNeubound3}
\begin{split}
\mathcal{P}^\Neu_\alpha(\lambda)&\ge \frac{\lambda^2}{4}-\lambda\left(\frac{1}{\pi}+\frac{2}{5}\cdot\frac{1}{2}\right)+\frac{2}{5}\left(\frac{1}{2}G_\lambda^{-1}\left(\frac{1}{4}\right)\right)+\frac{1}{2}G_\lambda^{-1}\left(\frac{1}{4}\right)-\frac{1}{4}\\
&\ge \frac{\lambda^2}{4}-\lambda\left(\frac{1}{\pi}+\frac{1}{5}\right)+\frac{7}{10}G_\lambda^{-1}\left(\frac{1}{4}\right)-\frac{1}{4}.
\end{split}
\end{equation}

We further recall \cite[Lemma 4.8]{FLPS} which after the substitution $c=\cos\sigma$ states that 
\begin{equation}\label{eq:Ginvb0}
G_\lambda^{-1}\left(\frac{1}{4}\right)\ge c\lambda
\end{equation} 
with $c\in(0,1]$ if 
\[
\lambda\ge q_1(c):=\frac{\pi}{4(\sqrt{1-c^2}-c\arccos c)}.
\]
Applying \eqref{eq:Ginvb0} in the right-hand side of \eqref{eq:PNeubound3}, we obtain 
\[
\mathcal{P}^\Neu_\alpha(\lambda) \ge \frac{\lambda^2}{4}+\left(\frac{7}{10}c-\frac{1}{\pi}-\frac{1}{5}\right)\lambda-\frac{1}{4},
\]
which is in turn greater than $\frac{\lambda^2}{4}$ assuming that
\[
\lambda>q_2(c):=\frac{5}{2\left(7c-2-\frac{10}{\pi}\right)}>0.
\]
We have therefore proved \eqref{eq:polyaneu} for all 
\begin{equation}\label{eq:implicitlambdaN}
\lambda>\min_{c\in\left(\frac{2\pi+10}{7\pi},1\right)}\max\{q_1(c),q_2(c),3\}.
\end{equation}

It is easily checked that the equation $q_1(c)=q_2(c)$ has a single root $c^*\in\left(\frac{2\pi+10}{7\pi},1\right)$, with $c^*\approx 0.78397$ and $q^*:=q_1(c^*)=q_2(c^*)\approx 8.2047$, see Figure \ref{fig:qs}.

\begin{figure}[htpb]
\centering
\includegraphics{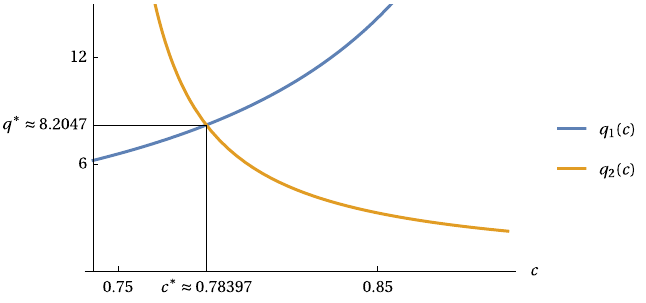}
\caption{The functions $q_1(c)$ and $q_2(c)$.\label{fig:qs}}
\end{figure}

To finish the proof of  Proposition \ref{prop:largelambda}, we drop the minimisation in the right-hand side of \eqref{eq:implicitlambdaN} and choose there the particular value $c=\frac{783}{1000}$. The method of verified rational approximations from \cite[\S 8]{FLPS} ensures that 
\[
q_2\left(\frac{783}{1000}\right)=\frac{2500\pi}{3481\pi-10000}>q_1\left(\frac{783}{1000}\right)>3.
\]
\end{proof}

\begin{figure}[htpb]
\centering
\includegraphics{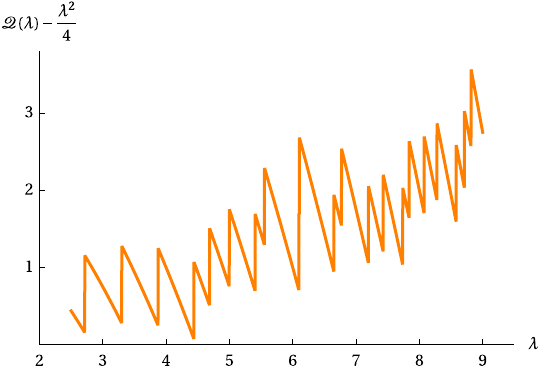}
\caption{The difference $\mathcal{Q}(\lambda)-\frac{\lambda^2}{4}$ plotted as a function of $\lambda\in\left[\frac{5}{2}, 9\right]$. The plot has been produced using floating-point arithmetic and is included for illustration only: it does not constitute a part of the proof. 
\label{fig:Qlambda}}
\end{figure}

The combination of Propositions \ref{prop:smallambda} and \ref{prop:largelambda} proves Theorem \ref{thm:polyan} except for $\lambda\in\left[2\sqrt{2},\frac{2500\pi}{3481\pi-10000}\right]\subset \left[\frac{5}{2},9\right]$.
In order to close the gap, we first obtain another lower bound for $\mathcal{P}^\Neu_\alpha(\lambda)$ by moving all arguments of $G_\lambda$ appearing in its definition to the right to the nearest integer or half-integer:
\[
\mathcal{P}^\Neu_\alpha(\lambda)>\sum_{m=0}^{\ceiling{\lambda}-1} \left(\entire{G_\lambda\left(m+\frac{1}{2}\right)+\frac{3}{4}}+\entire{G_\lambda\left(m+1\right)+\frac{3}{4}}\right)=:\mathcal{Q}(\lambda).
\]
Note that $\mathcal{Q}(\lambda)$ does not depend on $\alpha$. We have
\begin{prop}\label{prop:gap} For any $\lambda\in \left[\frac{5}{2},9\right]$,
\begin{equation}\label{eq:Qineq}
\mathcal{Q}(\lambda)-\frac{\lambda^2}{4}>0,
\end{equation}
see Figure \ref{fig:Qlambda} for an illustration.
\end{prop}

Proposition \ref{prop:gap} is proved using the same strategy of applying a  rigorous computer-assisted algorithm for estimating $\mathcal{Q}(\lambda)$ (which requires a finite number of steps and uses only integer arithmetic) as was done in the case of the Neumann Laplacian for the count  $\mathcal{P}^\Neu_0(\lambda)$ in \cite[\S 8]{FLPS}.  Let, for $x\in\mathbb{R}$, $\underline{x}\le x\le \overline{x}$, where $\underline{x}, \overline{x}\in\mathbb{Q}$ are some lower and upper rational approximations of $x$. For $\lambda\in\mathbb{Q}$ we compute a verified lower rational approximation $\underline{\mathcal{Q}}(\lambda)$ using the techniques of \cite[\S 8]{FLPS}. We now use the analogue of \cite[Lemma 8.1]{FLPS}: if \eqref{eq:Qineq} holds for some $\lambda_0\in\mathbb{Q}$, or, more precisely, if 
\[
\underline{e}\left(\lambda_0\right):=\underline{\mathcal{Q}}\left(\lambda_0\right)-\frac{\lambda_0^2}{4}>0,
\]
then \eqref{eq:Qineq} also holds for all $\lambda\in\left(\lambda_0+\underline{\delta}\left(\lambda_0\right)\right)$, where 
\[
\delta\left(\lambda_0\right):=\sqrt{\lambda_0^2+4\underline{e}\left(\lambda_0\right)}-\lambda_0.
\]
The results of the calculations are shown in Table \ref{table:1}, and the script we used is available at

\centerline{\url{https://michaellevitin.net/polya.html\#AB}.}

\begin{table}[htb]
{\setlength{\tabcolsep}{10pt}
\centering
\begin{longtable}{@{}>{$}c<{$}*{3}{>{$}c<{$}}@{}}
\toprule
\text{Step}&\lambda&\underline{e}(\lambda)&\underline{\delta}(\lambda)\\\addlinespace[1.5ex]
\midrule
1 & \frac{5}{2} & \frac{7}{16} & \frac{15}{46} \\\addlinespace[1.5ex]
2 & \frac{65}{23} & \frac{2123}{2116} & \frac{190}{299} \\\addlinespace[1.5ex]
 3 & \frac{45}{13} & \frac{679}{676} & \frac{7}{13} \\\addlinespace[1.5ex]
 4 & 4 & 1 & \frac{8}{17} \\\addlinespace[1.5ex]
 5 & \frac{76}{17} & \frac{290}{289} & \frac{210}{493} \\\addlinespace[1.5ex]
 6 & \frac{142}{29} & \frac{846}{841} & \frac{354}{899} \\\addlinespace[1.5ex]
 7 & \frac{164}{31} & \frac{964}{961} & \frac{328}{899} \\\addlinespace[1.5ex]
 8 & \frac{164}{29} & \frac{1686}{841} & \frac{600}{899} \\\addlinespace[1.5ex]
 9 & \frac{196}{31} & \frac{1928}{961} & \frac{505}{837} \\\addlinespace[1.5ex]
 10 & \frac{187}{27} & \frac{5855}{2916} & \frac{5}{9} \\\addlinespace[1.5ex]
 11 & \frac{202}{27} & \frac{1463}{729} & \frac{14}{27} \\\addlinespace[1.5ex]
 12 & 8 & 2 & \frac{14}{29} \\\addlinespace[1.5ex]
 13 & \frac{246}{29} & \frac{1691}{841} & \frac{467}{1015} \\\addlinespace[1.5ex]
 14 & \frac{313}{35} & \frac{14731}{4900} & \frac{883}{1365} \\\addlinespace[1.5ex]
 15 & \frac{374}{39}>9 &  &  \\\addlinespace[1.5ex]
  \bottomrule\\\addlinespace[2ex]
   \caption{Detailed output of the computer-assisted algorithm.}\label{table:1}
\end{longtable}}
\end{table}

\end{document}